\documentclass[12pt]{article}
\usepackage{amsmath,amssymb,amsthm,comment,cite}

\newtheorem{theorem}{Theorem}[section]

\def\barr{\begin{array}}
\def\earr{\end{array}}

\title{Finite groups with many isolated subgroups}
\author{Marius T\u arn\u auceanu}
\date{November 28, 2021}

\begin{document}

\maketitle

\begin{abstract}
Given a finite group $G$, we denote by $L(G)$ the subgroup lattice of $G$ and by ${\rm Isolated}(G)$ the set of isolated subgroups of $G$.
In this note, we describe finite groups $G$ such that $|{\rm Isolated}(G)|=|L(G)|-k$, where $k=0,1,2$.
\end{abstract}
\smallskip

{\small
\noindent
{\bf MSC 2020\,:} Primary 20D30; Secondary 20D60.

\noindent
{\bf Key words\,:} finite group, subgroup lattice, isolated subgroup, ${\rm CP}_1$-group.}

\section{Introduction}

Let $G$ be a finite group and $L(G)$ be the subgroup lattice of $G$. We say that a subgroup $H\in L(G)$ is \textit{isolated} in $G$ if for every $x\in G$ we have either $x\in H$ or $\langle x\rangle\cap H=1$. Groups with isolated subgroups were studied in \cite{1,2}. However, this concept appears much earlier (see for instance Section 66 of \cite{5} and the entry "isolated subgroup" in the Encyclopedia of Mathematics, c.f. \cite{11}).

Recently, the set ${\rm Isolated}(G)$ of isolated subgroups of $G$ has been investigated for certain classes of non-abelian $p$-groups (see \cite{4}) and for abelian groups (see \cite{10}). Also, the structure of finite $p$-groups with isolated minimal non-abelian subgroups and of finite $p$-groups with an isolated metacyclic subgroup has been determined in \cite{8}.

The current note deals with the following natural question:

\begin{center}
\textit{How large can the set ${\rm Isolated}(G)$ \!be}?
\end{center}

Our main result is stated as follows. We recall that a group is called a \textit{${\rm CP}_1$-group} if all its elements have prime order; the structure of such a group is given by \cite{3} (see also \cite{6}).

\begin{theorem}
Let $G$ be a finite group. Then
\begin{itemize}
\item[{\rm a)}] ${\rm Isolated}(G)=L(G)$ if and only if $G$ is a ${\rm CP}_1$-group;
\item[{\rm b)}] $|{\rm Isolated}(G)|=|L(G)|-1$ if and only if $G\cong\mathbb{Z}_{p^2}$ with $p$ a prime;
\item[{\rm c)}] $|{\rm Isolated}(G)|=|L(G)|-2$ if and only if $G\cong\mathbb{Z}_{p^3}$ with $p$ a prime or $G\cong\mathbb{Z}_{pq}$ with $p$ and $q$ distinct primes.
\end{itemize}
\end{theorem}\smallskip

The classification of finite groups $G$ such that
\begin{equation}
|{\rm Isolated}(G)|=|L(G)|-k\nonumber
\end{equation}can be continued for $k\geq 3$. We remark that cyclic groups of order $p^{k+1}$, where $p$ is a prime, or of order $p_1p_2\cdots p_k$, where $p_1,p_2,\dots,p_k$ are distinct primes, satisfy this property. In other words, for every positive integer $k$ there is a finite group with exactly $k$ non-isolated subgroups.\smallskip

For the proof of the above theorem, we need the following well-known results:

\begin{itemize}
\item[$\bullet$] A finite group with one or two maximal subgroups is cyclic.
\item[$\bullet$] A finite group has a unique minimal subgroup if and only if it is either a cyclic $p$-group or a generalized quaternion $2$-group.
\item[$\bullet$] A finite group has exactly two minimal subgroups if and only if it is either cyclic of order $p^mq^n$ with $p,q$ distinct primes or a direct product of a cyclic $p$-group of odd order and a generalized quaternion $2$-group.
\end{itemize}

By a \textit{generalized quaternion $2$-group} we mean a group of order $2^n$ for some positive integer $n\geq 3$, defined by
\begin{equation}
Q_{2^n}=\langle a,b \mid a^{2^{n-2}}= b^2, a^{2^{n-1}}=1, b^{-1}ab=a^{-1}\rangle.\nonumber
\end{equation}

Finally, we note that the dual problem of finding finite groups with few isolated subgroups seems to be also very interesting. We say that a finite group is \textit{isolated-simple} if it has no proper isolated subgroups. Some obvious examples of such groups are cyclic groups, generalized quaternion $2$-groups or abelian $p$-groups of type $\mathbb{Z}_{p^{\alpha_1}}\times\mathbb{Z}_{p^{\alpha_2}}\times\cdots\times\mathbb{Z}_{p^{\alpha_k}}$ with $2\leq\alpha_1\leq\alpha_2\leq...\leq\alpha_k$ (see Lemma 2.3 of \cite{10}).

\medskip\noindent{\bf Open problem.} Determine the structure of isolated-simple groups.
\medskip

Most of our notation is standard and will not be repeated here. Basic definitions and results on group theory can be found in \cite{5,9}. For subgroup lattice concepts we refer the reader to \cite{7}.\newpage

\section{Proof of Theorem 1.1}

a) Assume that ${\rm Isolated}(G)=L(G)$, but $G$ is not a ${\rm CP}_1$-group. Then there exists $a\in G$ with $o(a)=pq$, where $p$ and $q$ are possibly equal primes. Since $a\notin\langle a^p\rangle$ and $\langle a\rangle\cap\langle a^p\rangle=\langle a^p\rangle\neq 1$, it follows that $\langle a^p\rangle$ is not isolated, a contradiction.

Conversely, let $H\in L(G)$ and $a\notin H$. Then $\langle a\rangle\cap H\leq\langle a\rangle$ and $\langle a\rangle\cap H\neq\langle a\rangle$. Since $\langle a\rangle$ is of prime order, we get $\langle a\rangle\cap H=1$. This shows that $H\in{\rm Isolated}(G)$, as desired.
\medskip

b) Assume that ${\rm Isolated}(G)=L(G)\setminus\{H_0\}$. Then all conjugates of $H_0$ are also not isolated and therefore $H_0$ is normal in $G$.

Let $a\in G\setminus H_0$ such that $\langle a\rangle\cap H_0\neq 1$. Then $a\notin\langle a\rangle\cap H_0$ and 
\begin{equation}
\langle a\rangle\cap(\langle a\rangle\cap H_0)=\langle a\rangle\cap H_0\neq 1,\nonumber
\end{equation}which implies that $\langle a\rangle\cap H_0$ is not isolated in $G$. It follows that $\langle a\rangle\cap H_0=H_0$, i.e. $H_0\subseteq\langle a\rangle$. Clearly, $H_0$ cannot have proper subgroups because they would not be isolated, contradicting our assumption. Thus, $|H_0|=p$ with $p$ a prime. We also observe that any proper subgroup of $\langle a\rangle$ must coincide with $H_0$ because it is not isolated. This leads to $|\langle a\rangle|=p^2$.

Suppose that $G$ possesses a subgroup $H$ of order $p^3$ containing $\langle a\rangle$. Then $H$ is of one of the following types: $\mathbb{Z}_{p^3}$, $\mathbb{Z}_p\times\mathbb{Z}_{p^2}$, the non-abelian group of order $p^3$ and exponent $p^2$, $D_8$ or $Q_8$. It is easy to see that if $H\cong\mathbb{Z}_{p^3}$ then $\langle a\rangle$ is not isolated, while if $H\cong D_8$ then the non-cyclic subgroups of order $4$ of $H$ are not isolated, contrary to our assumption. In the other cases there exists $\langle b\rangle\leq H$ with $|\langle b\rangle|=p^2$ and $\langle b\rangle\cap\langle a\rangle=H_0$, proving that $\langle a\rangle$ is again not isolated, a contradiction. Thus, $\langle a\rangle$ is a Sylow $p$-subgroup of $G$.

Suppose now that $|G|$ has a prime divisor $q\neq p$ and let $K\leq G$ with $|K|=q$. Then $H_0K\leq G$, $a\notin H_0K$ and 
\begin{equation}
\langle a\rangle\cap H_0K=H_0\neq 1,\nonumber
\end{equation}showing that $H_0K$ is not isolated in $G$, a contradiction.

Hence $G=\langle a\rangle\cong\mathbb{Z}_{p^2}$, as desired.
\medskip

c) Assume that ${\rm Isolated}(G)=L(G)\setminus\{H_1,H_2\}$. Similarly, $H_1$ and $H_2$ are normal in $G$.

Let $a\in G\setminus H_1$ such that $\langle a\rangle\cap H_1\neq 1$. Then $\langle a\rangle\cap H_1$ is not isolated, implying that $\langle a\rangle\cap H_1\in\{H_1,H_2\}$. We distinguish the following two cases:\newpage

\medskip\noindent\hspace{15mm}{\bf Case 1.} $\langle a\rangle\cap H_1=H_1$\medskip

Then $H_1\subset\langle a\rangle$. Since proper subgroups of $\langle a\rangle$ are not isolated, we infer that they are at least one, namely $H_1$, and at most two, namely $H_1$ and $H_2$. So, we have two possibilities:

\medskip\noindent\hspace{20mm}{\bf Subcase 1.1.} $H_1$ is the unique proper subgroup of $\langle a\rangle$\medskip

Then $|\langle a\rangle|=p^2$ and $|H_1|=p$, where $p$ is a prime.

Suppose first that $H_1$ is the unique minimal subgroup of $G$. This implies that $G$ is either a cyclic $p$-group or a generalized quaternion $2$-group. Then it is easy to see that $G$ has exactly two non-isolated subgroups if and only if $G\cong\mathbb{Z}_{p^3}$.

Suppose now that $G$ possesses minimal subgroups different from $H_1$ and let $K$ be such a subgroup. Then $a\notin H_1K$ and 
\begin{equation}
\langle a\rangle\cap H_1K=H_1\neq 1,\nonumber 
\end{equation}that is $H_1K$ is not isolated in $G$. Since $H_1K\neq H_1$, we get $H_1K=H_2$. Thus $H_2$ contains all minimal subgroups of $G$.

Observe also that $\langle a\rangle$ is normal in $G$. Indeed, for every $x\in G$, if $\langle a\rangle^x\neq\langle a\rangle$ then $\langle a\rangle\cap\langle a\rangle^x=1$ because $\langle a\rangle^x$ is isolated in $G$. It follows that $\langle a^p\rangle^x$ is not isolated in $G$ and since $\langle a^p\rangle^x\neq H_1,H_2$ we get a contradiction.

Let $H\neq\langle a\rangle, H_2$ be a maximal subgroup of $\langle a\rangle H_2$. Then $H$ is isolated in $G$ and therefore $\langle a\rangle\cap H=1$. Since $\langle a\rangle H=\langle a\rangle H_2$, it follows that $H$ is a minimal subgroup of $G$, which leads to $H\subset H_2$. Thus $\langle a\rangle H_2$ has exactly two maximal subgroups, namely $\langle a\rangle$ and $H_2$, and consequently it is cyclic. This implies that $\langle a\rangle$ is not isolated in $G$, a contradiction.

\medskip\noindent\hspace{20mm}{\bf Subcase 1.2.} $H_1$ and $H_2$ are the unique proper subgroups of $\langle a\rangle$\medskip

Then we have either
\begin{equation}
|\langle a\rangle|=p^3 \mbox{ and } |H_1|=p, |H_2|=p^2, \mbox{ where } p \mbox{ is a prime}\nonumber
\end{equation}or
\begin{equation}
|\langle a\rangle|=pq \mbox{ and } |H_1|=p, |H_2|=q, \mbox{ where } p,q \mbox{ are distinct primes}.\nonumber
\end{equation}In both cases, if $G$ would have a minimal subgroup $K$ different from $H_1$ and $H_1,H_2$, respectively, then $a\notin H_1K$ and
\begin{equation}
\langle a\rangle\cap H_1K=H_1\neq 1,\nonumber 
\end{equation}i.e. $H_1K$ is non-isolated in $G$. Since obviously $H_1K\neq H_1,H_2$, this contradicts our hypothesis. Thus $H_1$ and $H_1,H_2$, respectively, are the unique minimal subgroups of $G$, proving that $G$ is one of the following groups:
\begin{itemize}
\item[-] $\mathbb{Z}_{p^m}$, where $p$ is a prime and $m\in\mathbb{N}$;
\item[-] $Q_{2^n}$, where $n\in\mathbb{N}$, $n\geq 3$;
\item[-] $\mathbb{Z}_{p^mq^n}$, where $p,q$ are distinct primes and $m,n\in\mathbb{N}$;
\item[-] $\mathbb{Z}_{p^m}\times Q_{2^n}$, where $p$ is an odd prime and $m,n\in\mathbb{N}$, $n\geq 3$.
\end{itemize}In all cases one can easily check that $G$ has exactly two non-isolated subgroups if and only if it coincides with $\langle a\rangle$, as desired.

\medskip\noindent\hspace{15mm}{\bf Case 2.} $\langle a\rangle\cap H_1=H_2$\medskip

Then $H_2\subseteq\langle a\rangle$ and $H_2\subset H_1$. If $H_2=\langle a\rangle$, then we get $\langle a\rangle\subseteq H_1$, a contradiction. Thus $H_2\neq\langle a\rangle$. Again, proper subgroups of $\langle a\rangle$ are not isolated in $G$. If $H_2$ is the unique proper subgroup of $\langle a\rangle$, then we proceed similarly with Subcase 1.1, while if $H_1$ and $H_2$ are the unique proper subgroups of $\langle a\rangle$, then we proceed similarly with Subcase 1.2.\medskip

This completes the proof.\qed

\vspace*{3ex}\small

\hfill
\begin{minipage}[t]{5cm}
Marius T\u arn\u auceanu \\
Faculty of  Mathematics \\
``Al.I. Cuza'' University \\
Ia\c si, Romania \\
e-mail: {\tt tarnauc@uaic.ro}
\end{minipage}

\end{document}